\documentclass[12pt,a4paper]{amsart}
\usepackage{a4wide}
\usepackage{amsfonts,amsthm,amsmath,amssymb}
\usepackage{graphicx}

\theoremstyle{plain}
\newtheorem{lem}{Lemma}[section]
\newtheorem{prop}[lem]{Proposition}
\newtheorem{thm}[lem]{Theorem}
\newtheorem*{thm*}{Theorem}
\newtheorem*{thmA}{Theorem A}
\newtheorem*{thmB}{Theorem B}
\newtheorem{cor}[lem]{Corollary}
\newtheorem*{cor*}{Corollary}

\newtheorem*{koebe}{The Koebe Distortion Theorem}
\newtheorem*{koebe-sph}{The Spherical Koebe Distortion Theorem}
\newtheorem*{cover}{Spherical Distortion Theorem for logarithmic tracts}
\newtheorem*{bowen-S}{Bowen's formula for meromorphic maps in $\SSS$}
\newtheorem*{bowen-B-hyp}{Bowen's formula for hyperbolic meromorphic maps in $\B$}
\newtheorem*{bowen-B}{Bowen's formula for meromorphic maps in $\B$}

\theoremstyle{definition}
\newtheorem{defn}[lem]{Definition}
\newtheorem*{defn*}{Definition}

\newtheorem*{ex*}{Example}
\newtheorem{rem}[lem]{Remark}
\newtheorem*{rem*}{Remark}

\theoremstyle{remark}

\DeclareMathOperator{\diam}{diam}
\DeclareMathOperator{\Area}{Area}
\DeclareMathOperator{\Arg}{Arg}
\DeclareMathOperator{\Sing}{Sing}

\newcommand{\C}{\mathbb C}
\newcommand{\D}{\mathbb D}
\newcommand{\clC}{\overline \C}
\newcommand{\clD}{\overline \D}
\newcommand{\R}{\mathbb R}
\newcommand{\Z}{\mathbb Z}

\newcommand{\F}{\mathcal F}
\newcommand{\B}{\mathcal B}
\newcommand{\DD}{\mathcal D}
\newcommand{\SSS}{\mathcal S}
\newcommand{\PP}{\mathcal P}

\begin{document}

\title{Bowen's formula for meromorphic functions}

\date{}

\author{Krzysztof Bara\'nski}
\address{Institute of Mathematics, University of Warsaw,
ul.~Banacha~2, 02-097 Warszawa, Poland}
\email{baranski@mimuw.edu.pl}

\author{Bogus{\l}awa Karpi\'nska}
\address{Faculty of Mathematics and Information Science, Warsaw
University of Technology, Pl.~Politechniki~1, 00-661 Warszawa, Poland}
\email{bkarpin@mini.pw.edu.pl}

\author{Anna Zdunik}
\address{Institute of Mathematics, University of Warsaw,
ul.~Banacha~2, 02-097 Warszawa, Poland}
\email{A.Zdunik@mimuw.edu.pl}

\thanks{Research partially supported by Polish MNiSW Grants N N201 0234 33 and N N201 607940 and EU FP6 Marie Curie Programme RTN CODY. The second author is partially supported by Polish PW Grant 504G 1120 0011 000.}
\subjclass[2000]{Primary 37F10, 37F35. Secondary 28A80.}

\begin{abstract}
Let $f$ be an arbitrary transcendental entire or meromorphic function in the class $\SSS$ (i.e.~with finitely many singularities). We show that the topological pressure $P(f,t)$ for $t > 0$ can be defined as the common value of the pressures $P(f,t, z)$ for all $z \in \C$ up
to a set of Hausdorff dimension zero. Moreover, we prove that $P(f,t)$ equals the supremum of the pressures of $f|_X$ over all invariant 
hyperbolic subsets $X$ of the Julia set, and we prove Bowen's formula for
$f$, i.e.~we show that the Hausdorff dimension of the radial Julia set
of $f$ is equal to the infimum of the set of $t$, for which $P(f,t)$ is
non-positive. Similar results hold for (non-exceptional) transcendental entire or meromorphic functions $f$ in the class $\B$ (i.e.~with bounded set of singularities), for which the closure of the post-singular set does not contain the Julia set. 
\end{abstract}

\maketitle

\section{Introduction}\label{sec:intro}

The thermodynamical formalism, developed by D.~Ruelle, R.~Bowen and P.~Walters
in the 1970's (see e.g.~\cite{Ru}) has provided a number of useful tools to
study the geometry and ergodic properties of invariant hyperbolic
subsets (conformal repellers) of the Julia set $J(f)$ of a rational map $f$ on the Riemann sphere. Recall that in this setting a
conformal repeller is a compact set $X \subset J(f)$, such that $X$ is
$f$-invariant (i.e.~$f(X) \subset X$) and $|(f^k)'|_X > 1$ for some $k
> 0$.  In particular, the celebrated Bowen formula (see \cite{Bo}) asserts 
that the Hausdorff dimension of a conformal repeller $X$ (e.g.~the Julia
set for a hyperbolic rational map) is equal to the unique zero 
of the pressure function $t \mapsto P(f|_X, t)$, where 
\[
P(f|_X, t) =
\lim_{n\to\infty} \frac{1}{n} \ln \sum_{w \in f^{-n}(z) \cap X}
|(f^n)'(w)|^{-t}
\]
for $z \in X$ is the topological pressure of $f|_X$ for the potential $\varphi = -t\ln|f'|$. 

A rational map $f$ is hyperbolic, if the closure (in $\clC$) of
the post-critical set (i.e.~the union of forward trajectories of
the critical values of $f$) is disjoint from the Julia set
of $f$. For transcendental meromorphic maps the definition is slightly different --- in this case we call $f$ hyperbolic, if the closure (in $\clC$) of the post-singular set $\PP(f)$ (i.e.~the union of forward trajectories of
the singular values of $f$) is disjoint from $J(f) \cup \{\infty\}$. Recall that the  set of singular
values of $f$, denoted by $\Sing(f)$, is the set of all finite singularities of
$f^{-1}$ (critical and asymptotic values).

In recent years, there have been more and
more attempts to generalise at least some results of the thermodynamical
formalism theory to the case of transcendental meromorphic
maps. However, this encounters some difficulties, due to lack of compactness,
infinite degree of the map and more complicated geometry. In
particular, the pressure function for the potential given above is
usually infinite.  

The first idea to overcome this obstruction was to consider hyperbolic
transcendental meromorphic maps $f$ which are periodic with some
period $T \in \C$. Then one can project $f$ to the cylinder $\C/T\Z$,
which (in some cases) makes the pressure function finite.

Using these ideas, K.~Bara\'nski in \cite{B} developed some elements of the
thermodynamical formalism (in particular Bowen's formula for the dimension of the Julia set) for
certain hyperbolic meromorphic maps of the form $f(z) = h(\exp(az))$,
where $a\in\C$ and $h$ is a rational function, in particular for the hyperbolic maps from the tangent family $\lambda \tan (z)$, $\lambda \in \C$. The results were then generalised by J.~Kotus and M.~Urba\'nski in \cite{KU} to the case of so-called regular Walters expanding conformal maps.

In \cite{UZ, UZ2}, M.~Urba\'nski and A.~Zdunik
created the thermodynamical formalism theory for hyperbolic maps in the
exponential family $f(z) = \lambda \exp(z), \;\lambda \in \C$. In
particular, they discovered that for these maps Bowen's formula has a
different form. More precisely, the unique zero
of the pressure function is equal not to the Hausdorff dimension of
the Julia set $J(f)$ (which is equal to $2$ for all parameters $\lambda$, as
proved by C.~McMullen in \cite{McM}), but to the Hausdorff dimension
of the radial Julia set $J_r(f)$. The set $J_r(f)$ is, by
definition, the set of such $z\in J(f)$ for which there exists
$r=r(z)>0$ and a sequence $n_j\to\infty$, such that a holomorphic 
branch of $f^{-n_j}$ sending $f^{n_j}(z)$ to $z$ is well-defined on
the disc with respect to the spherical metric in $\clC$, centred at
$f^{n_j}(z)$ of radius $r$. In \cite{UZ, UZ2} it was proved that for hyperbolic
exponential maps the Hausdorff dimension of $J_r(f)$ is greater than
$1$ and smaller than $2$ and varies analytically with respect to the parameter
$\lambda$. The estimate holds also for some non-hyperbolic exponential maps, as proved in \cite{UZ3}. This shows that the radial Julia set $J_r(f)$ can be essentially smaller than the whole Julia set $J(f)$. 
Note that in \cite{UZ, UZ2} the set $J_r(f)$ was
defined in a (formally) different way, as the set of points in $J(f)$,
which do not escape to $\infty$ under iterates of $f$. However, it is
easy to see that the two definitions are equivalent for hyperbolic
exponential maps.

In \cite{MU1,MU2}, V.~Mayer and M.~Urba\'nski developed the thermodynamical formalism theory for hyperbolic transcendental meromorphic maps of finite order  
with the so-called balanced derivative growth condition. This condition is satisfied e.g.~when  
\[
c^{-1} (1 + |z|^\alpha)(1 + |f(z)|^\beta) \leq |f'(z)| \leq c 
(1 + |z|^\alpha)(1 + |f(z)|^\beta)
\]
for $z \in J(f) \setminus f^{-1}(\infty)$ and some $\alpha \in \R$,
$c, \beta > 0$ with $\alpha + \beta > 0$. Then in the definition of
the pressure instead of the standard derivative of $f^n$ one considers the
derivative with respect to the metric $d\sigma = \frac{dz}{1 +
|z|^\beta}$. Also in this case the unique zero of the pressure function is equal to the Hausdorff dimension of $J_r(f)$.
Among examples of maps satisfying the balanced derivative growth
condition are hyperbolic functions of the form $f(z) = P(\exp(Q(z))$,
where $P, Q$ are polynomials and hyperbolic (co)sine, tangent and
elliptic functions.  
The approach of \cite{MU1,MU2,UZ,UZ2} was based on a construction of a suitable conformal measure and an absolutely continuous invariant measure. Such a construction requires additional restrictive assumptions (such as balanced growth, finite order and hyperbolicity).

In this paper, using another approach, we show that Bowen's formula in its new form is actually
satisfied for all transcendental meromorphic maps in the class $\SSS$ and for a wide class of maps from the class $\B$.  
In particular, we need no additional restrictive conditions (like balanced growth or finite order). What is more, our proof works for non-hyperbolic maps as well.
Recall that the Speiser class $\SSS$ consists of transcendental meromorphic maps for which the set of singular values $\Sing(f)$ is finite. The Eremenko-Lyubich class $\B$ consists of transcendental meromorphic functions for which  $\Sing(f)$ is bounded. For results concerning the dynamics of maps from the classes $\SSS$ and $\B$ refer e.g.~to \cite{Ber, EL}. 

In this paper we prove the following.

\begin{thmA} For every transcendental entire or meromorphic map $f$ in the class $\SSS$ and every $t > 0$ the topological pressure 
\[
P(f,t) = P(f,t,z)=\lim_{n\to\infty}\frac 1 n\ln\sum_{w\in
f^{-n}(z)}|(f^n)^*(w)|^{-t}
\]
$($where ${}^*$ denotes the derivative with respect to the spherical metric$)$ 
exists $($possibly equal to $+\infty)$ and is independent of $z \in \C$ up to an exceptional set of Hausdorff dimension zero $($consisting of points quickly approximated by the forward orbits of singular values of $f)$. We have
\[
P(f, t) = P_{hyp}(f, t),
\] 
where $P_{hyp}(f, t)$ is the supremum of the pressures $P(f|_X, t)$
over all transitive isolated conformal repellers $X \subset J(f)$. The
function $t \mapsto P(f, t)$ is non-increasing and convex when it is
finite and satisfies $P(f,2) \leq 0$.
The following version of Bowen's formula holds:
\[
\dim_H J_r(f)  = \dim_{hyp} J(f) = \delta(f),
\]
where $\delta(f) = \inf \{t > 0: P(t) \leq 0\}$.
\end{thmA}

A conformal repeller $X$ is transitive, if for all non-empty sets $U, V$ open in $X$ we have $f^n(U) \cap V \neq \emptyset$ for some $n \geq 0$; $X$ is isolated, if there exists a neighbourhood $W$ of $X$, such that for every $z \in W\setminus X$ there exists $n >0$ with $f^n(z) \notin W$.

The hyperbolic dimension of the Julia set $J(f)$ (denoted
$\dim_{hyp}$), is defined as the supremum of the Hausdorff dimensions
(denoted $\dim_H$) of all conformal repellers contained in $J(f)$.
Obviously, the hyperbolic dimension is not greater than the Hausdorff
dimension. In \cite{R}, L.~Rempe proved that for transcendental
meromorphic maps, the hyperbolic dimension of the Julia set $J(f)$
coincides with the Hausdorff dimension of the radial limit set
$J_r(f)$. 

An analogue of Theorem~A holds for (non-exceptional) transcendental
meromorphic maps $f$ in the class $\B$, for which the closure of the
post-singular set $\PP(f)$ does not contain the whole Julia set, in particular for all hyperbolic maps from the class $\B$. 

We will call $f$ exceptional, if there exists a (Picard) exceptional
value $a$ of $f$, such that $a \in J(f)$ and $f$ has a non-logarithmic
singularity over $a$. Otherwise, we will say that $f$ is
non-exceptional. 

\begin{thmB} For every non-exceptional transcendental entire or meromorphic map
$f$ in the class $\B$, such that $J(f) \setminus \overline{\PP(f)} \neq \emptyset$ $($in particular, for every hyperbolic map in $\B)$ and every $t > 0$ the topological pressure 
\[
P(f,t) = P(f,t,z)=\lim_{n\to\infty}\frac 1 n\ln\sum_{w\in
f^{-n}(z)}|(f^n)^*(w)|^{-t}
\]
exists $($possibly equal to $+\infty)$ and is independent of
$z \in J(f) \setminus \overline{\PP(f)}$, which is an open dense
subset of $J(f)$. We have
\[
P(f, t) = P_{hyp}(f, t).
\] 
The function $t \mapsto P(f, t)$ is non-increasing and convex when it is finite
and satisfies $P(f,2) \leq 0$. Bowen's formula holds:
\[
\dim_H J_r(f)  = \dim_{hyp} J(f) = \delta(f).
\]
If, additionally, $f$ is hyperbolic, then $P(f, t) > 0$ for every $0 < t <\delta(f)$ and $P(f, t) < 0$ for every $t > \delta(f)$.
\end{thmB}

We were inspired by the papers by F.~Przytycki, J.~Rivera Letelier and
S.~Smirnov \cite{P, PRS}, where they developed the theory of the
pressure for arbitrary (not necessarily hyperbolic) rational
maps. In these papers they proved that the pressure $P(f, t)$ (with
the derivative of $f^n$ taken in the spherical metric) for such maps can be
defined as the common value of $P(f, t, z)$ for all $z \in \clC$ up to a
set of Hausdorff dimension zero and they showed the equivalence of
various kinds of pressures. In our paper, we prove that a similar theory can be developed in the case of transcendental meromorphic maps.
Some parts of our proofs repeat arguments and ideas used by
in \cite{P, PRS}. For completeness, we include these parts indicating suitable references.

Another source of inspiration was the paper \cite{S} by
G.~Stallard containing ideas which are very close to the notion of
the pressure for hyperbolic transcendental meromorphic maps in the 
class $\B$. 

The plan of the paper is as follows. After setting notation in
Section~\ref{sec:notat}, in Section~\ref{sec:lem} we prove a number of
technical facts. The most important one is the Spherical Distortion
Theorem for logarithmic tracts, which estimates the spherical
distortion of inverse branches of a holomorphic universal covering
of a punctured disc in the Riemann sphere. This theorem has some interest
in itself, since it provides useful estimates in a general
setup, e.g.~for inverse branches of a map $f \in \B$ with a finite
number of poles or, more generally, a map with a logarithmic tract (in the
sense of \cite{BRS}) near infinity. In Section~\ref{sec:pressure} we
define the pressure and introduce a notion of Good Pressure Starting
(GPS) points, i.e.~the points $z \in\C$ for which the pressure
$P(f,t,z)$ has good properties. Then in Section~\ref{sec:pressure-S} we
prove that for the maps in $\SSS$ the pressure $P(f,t,z)$ exists, is independent of the starting point $z$ within the set of GPS points and equals $P_{hyp}(f,t)$. In Section~\ref{sec:bowen-S} we
state Bowen's formula and complete the proof of Theorem~A. The last Section~\ref{sec:pressure-B} deals with the maps from class $\B$ and proves Theorem~B.

\section{Notation}\label{sec:notat}

Let $f$ be a transcendental meromorphic function (we treat entire functions as meromorphic). In what follows, we use the following notation: 
\begin{align*}
\Sing(f)&= \{z \in \C :z \text{ is a finite singularity of } f^{-1}\},\\
\PP_n(f) & = \bigcup_{k = 0}^{n-1} f^k(\Sing(f))= \Sing(f^n),\\
\PP(f)&=\bigcup_{k = 0}^\infty f^k(\Sing(f)).
\end{align*}
Here and in the sequel we use the symbol $f^k(A)$ to denote the image
under $f^k$ of the set of points in a set $A$ for which $f^k$ is
defined.  

We consider the standard spherical metric on the Riemann sphere $\clC$ defined by
\[
ds = \frac{2 \; dz}{1+|z|^2} 
\]
and denote by $d(z_1, z_2)$ the spherical distance between $z_1$ and $z_2$. 
We have 
\[
d(z_1, z_2) = 2 \arctan \left|\frac{z_1-z_2}{1+\bar z_1 z_2}\right|.
\]
We write $\D(z_0, r)$ (respectively $\DD(z_0, r)$) for the disc in the Euclidean (respectively spherical) metric, centred at $z_0$, of radius $r$, i.e.
\[
\D(z_0, r) = \{z\in\C: |z-z_0|<r\}, \qquad  \DD(z_0, r) = \{z\in\clC: d(z,z_0)<r\}. 
\]
For short, we write $\D(r) = \D(0, r)$. By $\clD(z_0,r)$ we denote the closed disc. For a set $A \subset \clC$ we write $d(z, A) = \inf\{d(z, w): w \in A\}$ and $\DD(A, r) = \{z\in\clC: d(z, A) < r\}$.
The diameter taken with respect to the spherical metric will be denoted by $\diam_{sph}$. The derivative of a holomorphic map $g$ with respect to the spherical metric will be denoted by $g^*$, while the standard derivative is denoted by $g'$. By definition,
\[
|g^*(z)| = \frac{(1+|z|^2)|g'(z)|}{1+|g(z)|^2}.
\]

\section{Distortion on logarithmic tracts}\label{sec:lem}

In this section we prove some technical lemmas that will be used in the proofs of the main results. Note that the constants $c, c_1, c_2$ etc. appearing in the lemmas may have different meanings.

Recall first the classical Koebe Distortion Theorem (see e.g.~\cite{CG}).

\begin{koebe}
Let $g: \D(z_0,r) \to \C$ be a univalent holomorphic map, for some $z_0 \in\C$ and $r > 0$. Then for every $z \in \D(z_0,r)$, if $|z - z_0| \leq \lambda r$ for some $0 < \lambda < 1$, then 
\begin{gather*}
\frac{1 - \lambda}{(1 + \lambda)^3} \leq \frac{|g'(z)|}{|g'(z_0)|} \leq \frac{1 + \lambda}{(1 - \lambda)^3},\\
\frac{\lambda r}{(1 + \lambda)^2}|g'(z_0)| \leq |g(z) - g(z_0)| \leq \frac{\lambda r}{(1 - \lambda)^2}|g'(z_0)|.
\end{gather*}
\end{koebe}

Now we prove an analogue of the Koebe Distortion Theorem for the spherical metric. This is a kind of folklore and it has already appeared (without proof) in several papers, but we think it is useful to present it here in a complete form.

\begin{koebe-sph}
Let $0 < r_1, r_2 < \diam_{sph} \clC$. Then there exists a constant $c>0$ depending only on $r_1, r_2$, such that for every spherical disc $D = \DD(z_0, r)$ and every univalent holomorphic map $g: D \to \clC$ with $z_0 \in \clC$, $\diam_{sph} D < r_1$ and $\diam_{sph}(\clC\setminus g(D))>r_2$, if $z_1, z_2 \in \DD(z_0, \lambda r)$ for some $0 < \lambda < 1$, then 
\[
\frac{|g^*(z_1)|}{|g^*(z_2)|} \leq\frac{c}{(1 - \lambda)^4}.
\]
\end{koebe-sph}

\begin{proof}
We denote by $c_1, c_2, \ldots$ constants depending only on $r_1, r_2$. Consider 
$z_1, z_2 \in \DD(z_0, \lambda r)$. Let $\varphi$ be an isometry of $\clC$ with respect to the spherical metric (i.e.~a rotation of the Riemann sphere), such that $\varphi(z_0) = 0$. Then $\varphi(D) = \DD(0,r) = \D(\varrho)$ and $\varphi(z_1), \varphi(z_2) \in \DD(0,\lambda r) = \D(\varrho_1)$ for some $\varrho_1 < \varrho < c_1$, such that 
\[
1 - \varrho_1/\varrho > c_2(1-\lambda).
\]
Since $\diam_{sph}(\clC\setminus g(D))>r_2$, there exists a point $w \in \clC\setminus g(D)$, such that $d(g(z_2), w) > r_2/2$. Take another spherical isometry $\psi$ with $\psi(w) = \infty$. Then $|\psi(g(z_2))| < c_3$. Define $\tilde g = \psi \circ g \circ \varphi^{-1}$ on $\D(\varrho)$. Then $\tilde g$ is a holomorphic univalent map into $\C$, such that $|\tilde g(\varphi(z_2))| < c_3$. By the classical Koebe Distortion Theorem for $\tilde g$,
\begin{multline*}
\frac{|g^*(z_1)|}{|g^*(z_2)|} = \frac{|\tilde g^*(\varphi(z_1))|}{|\tilde g^*(\varphi(z_2))|} =  \frac{(1+|\varphi(z_1)|^2)(1+|\psi(g(z_2))|^2)|\tilde g'(\varphi(z_1))|}{(1+|\varphi(z_2)|^2)(1+|\psi(g(z_1))|^2)|\tilde g'(\varphi(z_2))|} \\
< (1+c_1^2)(1+c_3^2) \frac{(1+\varrho_1/\varrho)^4}{(1-\varrho_1/\varrho)^4} < \frac{16 (1+c_1^2)(1+c_3^2)}{c_2^4(1-\lambda)^4},
\end{multline*}
which ends the proof.
\end{proof}

In \cite{BRS}, a general notion of a logarithmic tract (over infinity) was considered. In this paper we consider logarithmic tracts over any value $a \in \clC$. We recall the definition.

\begin{defn}
Suppose that $U \subset \C$ is an unbounded simply connected domain, such that the boundary of $U$ in $\C$ is a smooth open simple arc. Let $a \in \clC$ and $0 < r < \diam_{sph}\clC$. If $f:\overline{U} \to \C$ is a continuous map, holomorphic on $U$, such that $d(f(z), a) = r$ for every $z$ in the boundary of $U$ and $f$ on $U$ is a universal covering of $\DD(a, r) \setminus \{a\}$, then we call $U$ a logarithmic tract of $f$ over $a$.
\end{defn}

\begin{rem}\label{rem:log}
If $f$ is a meromorphic map on $\C$, then $f$ has a logarithmic tract over $a \in \clC$ if and only if $a$ is a logarithmic asymptotic value of $f$ (i.e.~$f$ has a logarithmic singularity over $a$). Note that if $f \in \B$ is entire or has a finite number of poles, then every component of $f^{-1}(V)$, where $V = \{z\in\C: |z| > R\}$ for sufficiently large $R$, is a logarithmic tract of $f$ over $\infty$.
\end{rem}

The definition of a logarithmic tract implies immediately the following.

\begin{lem}[\mbox{\cite[Lemma~6.1]{BRS}}]
Let $U$ be a logarithmic tract of $f$ over $a$, such that $f(U) = \DD(a, r) \setminus \{a\}$. Then for every $r' < r$, the set $U' = U \cap f^{-1}(\DD(a, r'))$ is also a logarithmic tract of $f$ and the boundary of $U'$ is an analytic open simple arc.
\end{lem}

This implies that, by diminishing $r$ if necessary, we can assume that $U$ does not contain $0$. Then, since $f$ is a universal covering, we can lift it to the logarithmic coordinates, i.e.~define a map 
\begin{equation}\label{eq:F}
F: \bigcup_{s \in \Z} (\log U + 2\pi i s) \to \{z\in\C:\textup{Re}(z) > \ln r\},
\end{equation}
where $\log$ is a branch of the logarithm on $U$, such that $\exp\circ F = f\circ \exp$. The map $F$ is periodic with period $2\pi i$ and maps  each set $\log U + 2\pi i s$ univalently onto the half-plane $\{z\in\C:\textup{Re}(z) > \ln r\}$.

The following estimation is well-known (see e.g.~\cite{BRS,EL}).

\begin{lem}\label{lem:EL}
Let $U$ be a logarithmic tract of $f: U \to V$ over $\infty$, where $V = \{z\in\C: |z| > R\}$ for some $R > 0$ and $0 \notin U$. Then for every $z \in U$ and every $w \in \bigcup_{s \in \Z} (\log U + 2\pi i s)$,
\[
|f'(z)| > \frac{|f(z)|(\ln|f(z)| -\ln R)}{4\pi|z|}, \qquad |F'(w)| > \frac{\textup{Re}(F(w)) - \ln R}{4\pi}.
\]
\end{lem}

Using Lemma~\ref{lem:EL} and the definition of the spherical metric we obtain the following corollary.

\begin{cor}\label{cor:tract_infty}
Let $R,L > 1$. Then there exists a constant $c > 0$ depending only on $R, L$, such that for every logarithmic tract $U \subset \C$ of $f:U \to V$ over $\infty$, where $V = \{z\in\C:|z| > R\}$ and $0 \notin U$, for every $z \in U$ with $|f(z)| > LR$ we have
\[
|f^*(z)| > c \frac{|z| \ln |f(z)|}{|f(z)|}.
\]
\end{cor}

In fact, we can reformulate the result in the general case of a tract over $a \in \clC$.

\begin{cor}\label{cor:tract}
Let $r,\lambda \in (0, 1)$. Then there exists a constant $c > 0$ depending only on $r, \lambda$, such that for every logarithmic tract $U \subset \C$ of $f:U \to V$ over $a \in \clC$, where $V = \DD(a,r)\setminus\{a\}$ and $0 \notin U$, for every $z \in U$ with $d(f(z), a) < \lambda r$ we have
\[
|f^*(z)| > c |z| d(f(z), a) \ln \frac{1}{d(f(z), a)}.
\]
\end{cor}
\begin{proof}It is sufficient to consider $h \circ f$, where $h$ is a spherical  isometry of $\clC$ such that $h(a) = \infty$, and notice that $(1/\tan(1/2))/|z| < d(z, \infty) < 2/|z|$ for $d(z,\infty) < 1$. 
\end{proof}

\begin{cor} \label{cor:|g|_infty}
Let $R, L > 1$. Then there exists a constant $c > 0$ depending only on $R, L$, such that for every logarithmic tract $U \subset \C$ of $f:U \to V$ over $\infty$, where $V = \{z\in\C:|z| > R\}$ and $0 \notin U$, for every $z_1, z_2 \in V$ with $|z_1| \geq |z_2| \geq LR$ and every branch $g$ of $f^{-1}$ in a neighbourhood of $z_1$ $($or $z_2)$, we have
\[
c^{-1}\left(\frac{\ln|z_1|}{\ln|z_2|}\right)^{-4\pi} < \frac{|g(z_1)|}{|g(z_2)|} < c\left(\frac{\ln|z_1|}{\ln|z_2|}\right)^{4\pi} 
\]
for some extension of the branch $g$ to a neighbourhood of $z_2$ $($or $z_1)$.
\end{cor}
\begin{proof} Assume that $g$ is defined in a neighbourhood of $z_1$ (the other case is symmetric). Let $w_1 = F(\log g(z_1))$ for the map $F$ from \eqref{eq:F}. By the definition of $F$, we have $w_1 = \ln |z_1| + 2 \pi i s + i \theta$ for some $s \in \Z$, $\theta \in \R$. Let $w_2 = \ln |z_2| + 2\pi i s + i \Arg(z_2)$, where $\Arg(z_2)$ is chosen such that $|\Arg(z_2) - \theta| < 2\pi$, and let $w_3 = \textup{Re}(w_1) + i\textup{Im}(w_2)$. Note that $\textup{Re}(w_1) = \textup{Re}(w_3) = \ln|z_1| \geq \ln(LR) > 0$ and $\textup{Re}(w_2) = \ln|z_2| \geq \ln(LR) > 0$.

Set $G = (F|_{\log U})^{-1}$. The branch $G$ is defined on the half-plane $\{z\in\C:\textup{Re}(z) >\ln R\}$, in particular in a neighbourhood of the curve $\gamma$, which is the union of two straight line segments connecting respectively $w_1$ to $w_3$ and $w_3$ to $w_2$. By Lemma~\ref{lem:EL}, $|G'(w)| < \frac{4\pi}{\textup{Re}(w)-\ln R}$ for $\textup{Re}(w) > \ln R$. Thus, integrating $G'$ along $\gamma$, we obtain
\begin{multline*}
|\textup{Re}(G(w_1)) - \textup{Re}(G(w_2))| \leq 
|G(w_1)-G(w_2)|\leq |G(w_1)-G(w_3)|+ |G(w_3)-G(w_2)| \\< \frac{8\pi^2}{\textup{Re}(w_1)-\ln R} + 4\pi \ln\frac{\textup{Re}(w_1)-\ln R}{\textup{Re}(w_2) - \ln R} \leq \frac{8\pi^2}{\ln L} + 4\pi \ln\frac{\ln (LR)}{\ln L}+ 4\pi \ln\frac{\textup{Re}(w_1)}{\textup{Re}(w_2)},
\end{multline*}
so
\begin{equation}\label{eq:|G|}
\textup{Re}(G(w_2)) - c_1 - 4\pi\ln\frac{\ln|z_1|}{\ln|z_2|} < \textup{Re}(G(w_1)) < \textup{Re}(G(w_2)) + c_1 + 4\pi \ln\frac{\ln|z_1|}{\ln|z_2|},
\end{equation}
where the constant $c_1 > 0$ depends only on $R,L$. 
Note that the map $\exp$ is univalent in a neighbourhood of the curve $\gamma$, because $|\textup{Im}(w_1) - \textup{Im}(w_2)| < 2\pi$. Hence, we can extend $g$ along the curve $\exp \gamma$ to a neighbourhood of $z_2$, such that $g \circ \exp = \exp \circ G$ on $\gamma$. Then $\textup{Re}(G(w_j)) = \ln |g(z_j)|$ for $j = 1, 2$. Substituting this into \eqref{eq:|G|}, we get the assertion. 
\end{proof}

As previously, we generalise the above corollary to the case of logarithmic tracts over an arbitrary point $a \in \clC$.

\begin{cor} \label{cor:|g|}
Let $r,\lambda \in (0, 1)$. Then there exists a constant $c > 0$ depending only on $r,\lambda$, such that for every logarithmic tract $U \subset \C$ of $f:U \to V$ over $a \in \clC$, where $V = \DD(a,r)\setminus\{a\}$ and $0 \notin U$, for every $z_1, z_2 \in U$ with $d(z_1,a) \leq d(z_2,a) \leq \lambda r$ and every branch $g$ of $f^{-1}$ in a neighbourhood of $z_1$ $($or $z_2)$, we have
\[
c^{-1}\left(\frac{\ln(1/d(z_1,a))}{\ln(1/d(z_2,a))}\right)^{-4\pi} < \frac{|g(z_1)|}{|g(z_2)|} < c\left(\frac{\ln(1/d(z_1,a))}{\ln(1/d(z_2,a))}\right)^{4\pi} 
\]
for some extension of the branch $g$ to a neighbourhood of $z_2$ $($or $z_1)$.
\end{cor}
\begin{proof}
The result follows in the same way as Corollary~\ref{cor:tract}.
\end{proof}

Now we prove a distortion theorem (in the spherical metric) for inverse branches of a map $f$ on a logarithmic tract. It is an improvement of \cite[Lemma~2.6]{S}.

\begin{cover} Let $r,\lambda \in (0, 1)$. Then there exists a constant $c > 0$ depending only on $r, \lambda$, such that for every logarithmic tract $U \subset \C$ of $f:U \to V$ over $a \in \clC$, where $V = \DD(a,r)\setminus\{a\}$ and $0\notin U$, for every $z_1, z_2 \in V$ with $d(z_1, a) \leq d(z_2, a) \leq \lambda r$ and every branch $g$ of $f^{-1}$ 
in a neighbourhood of $z_1$ $($or $z_2)$, we have
\begin{equation}\label{eq:sta}
c^{-1} \frac{d(z_2, a)}{d(z_1, a)} \left( \frac{\ln(1/d(z_1, a)}{\ln(1/d(z_2, a))}\right)^{-3} \leq \frac{|g^*(z_1)|}{|g^*(z_2)|} \leq c \frac{d(z_2, a)}{d(z_1, a)} \frac{\ln(1/d(z_1, a))}{\ln(1/d(z_2, a))}
\end{equation}
for some extension of the branch $g$ to a neighbourhood of $z_2$ $($or $z_1)$. In fact, the branch $g$ has two extensions $g_1$, $g_2$ such that for every $z_2$ $($or $z_1)$ as above, \eqref{eq:sta} holds for $g = g_1$ or $g = g_2$. 
\end{cover}

In particular, for $a = \infty$ we obtain the following.

\begin{cor}\label{cor:sta} Let $R, L > 1$. Then there exists a constant $c > 0$ depending only on $R, L$, such that for every logarithmic tract $U\subset \C$ of $f:U \to V$ over $\infty$, where $V = \{z\in\C: |z| > R\}$ and $0 \notin U$, for every $z_1, z_2 \in V$ with $|z_1| \geq |z_2| \geq LR$ and every branch $g$ of $f^{-1}$ in a neighbourhood of one of the points $z_1, z_2$, we have
\[
c^{-1} \frac{|z_1|}{|z_2|} \left(\frac{\ln |z_1|}{\ln |z_2|}\right)^{-3} \leq \frac{|g^*(z_1)|}{|g^*(z_2)|} \leq c \frac{|z_1|}{|z_2|} \frac{\ln |z_1|}{\ln |z_2|},
\]
for some extension of the branch $g$.
\end{cor}

\begin{proof}[Proof of Spherical Distortion Theorem for logarithmic tracts]

Note first that (as in the proof of Corollary~\ref{cor:tract}) taking $h \circ f$, where $h$ is a spherical isometry of $\clC$ such that $h(a) = \infty$, it is sufficient to consider the case $a = \infty$, i.e.~to prove Corollary~\ref{cor:sta}. 

To prove Corollary~\ref{cor:sta}, we first show that if additionally we assume $|z_1| = |z_2|$, then 
\begin{equation}\label{eq:first_step}
c_1^{-1} < \frac{|g^*(z_1)|}{|g^*(z_2)|} < c_1
\end{equation}
for some constant $c_1 > 0$ depending only on $R,L$. To show \eqref{eq:first_step}, note that by the definition of the spherical metric, we have 
\[
d(z_1, \infty) > \frac{c_2}{|z_1|}, \quad d(z_1, \clD(R)) > c_2, \quad d(z_1,z_1e^{it}) \leq c_2\frac{|1-e^{it}|}{|z_1|}
\]
for every $t \in \R$, where the constant $c_2 > 0$ depends only on $R,L$. This implies that 
\[
D = \DD(z_1, 2d(z_1, z_1e^{it})) \subset V \quad\text{ for } |t| \leq \delta,
\]
if $\delta$ is chosen sufficiently small (depending only on $R,L$). Extend the branch $g$ to $D$. Note that $\diam_{sph}D < \diam_{sph}\clC/2$ if $\delta$ is sufficiently small. Moreover, $\diam_{sph}(\clC\setminus g(D)) \geq \diam_{sph}(\clC\setminus U) = \diam_{sph}\clC$ since $0, \infty \notin U$. Hence, we can use the Spherical Koebe Distortion Theorem to obtain
\[
c_3^{-1} < \frac{|g^*(z_1e^{it})|}{|g^*(z_1)|} < c_3
\]
for every $|t| \leq \delta$, where the constant $c_3 > 0$ depends only on $R,L$. Since $z_2 = z_1 e^{i\theta}$ for some $\theta \in \R$ with $|\theta| \leq \pi$, we obtain \eqref{eq:first_step} with $c_1 = c_3^{\pi/\delta}$.

By \eqref{eq:first_step}, to prove Corollary~\ref{cor:sta} we can assume $z_1, z_2 \in \R^+$, such that $LR \leq z_2 \leq z_1$. Since $0 \notin U$, we can lift $f$ to the map $F$ from \eqref{eq:F}. 
By the definition of $F$, we have $F(\log g(z_1)) = w_1$, where $w_1 = \ln z_1 + 2 \pi i s$ for some $s \in \Z$. Let 
\[
H(w) = \exp G(w), \qquad \tilde H(w) = \frac{1}{H(w)} = \frac{1}{\exp G(w)},
\]
where $G = (F|_{\log U})^{-1}$. Then the maps $H$ and $\tilde H$ are well-defined and univalent in the half-plane $\{z\in\C:\textup{Re}(z) > \ln R\}$, in particular in the disc $\D(w_1, \ln (z_1/R))$. Hence, by the classical Koebe Distortion Theorem, 
\begin{equation}\label{eq:second_step}
\frac{1}{8} \frac{\ln( z_2/R)}{\ln( z_1/R)} \leq 
\frac{|H'(w_2)|}{|H'(w_1)|} \leq 2 \left(\frac{\ln( z_1/R)}{\ln( z_2/R)}\right)^3, \quad \frac{1}{8} \frac{\ln( z_2/R)}{\ln( z_1/R)} \leq
\frac{|\tilde H'(w_2)|}{|\tilde H'(w_1)|} \leq 2 \left(\frac{\ln( z_1/R)}{\ln( z_2/R)}\right)^3,
\end{equation}
where $w_2 = \ln z_2 + 2 \pi i s$.  

Note that in the neighbourhood of $z_1$, we have $g = H\circ \log_1$, where $\log_1$ is the branch of logarithm sending $z_1$ to $w_1$. Extending this branch to a neighbourhood of $(R,+\infty)$, we extend $g$ to a neighbourhood of $z_2$. Then $H = g\circ \exp$ in some neighbourhoods of $w_1$ and $w_2$, so $|H'(w_j)| = |g'(z_j)|z_j$ and $|\tilde H'(w_j)| = |g'(z_j)|z_j/|g(z_j)|^2$
for $j = 1, 2$. Hence, we can rewrite \eqref{eq:second_step} as
\begin{align}
\label{eq:1}
\frac{1}{2} \frac{z_2}{z_1}\left(\frac{\ln( z_2/R)}{\ln( z_1/R)}\right)^3 &\leq \frac{|g'(z_1)|}{|g'(z_2)|} \leq
8 \frac{z_2}{z_1} \frac{\ln( z_1/R)}{\ln( z_2/R)},\\
\label{eq:2}
\frac{1}{2} \frac{z_2}{z_1}\left(\frac{\ln( z_2/R)}{\ln( z_1/R)}\right)^3 \frac{|g(z_1)|^2}{|g(z_2)|^2}&\leq \frac{|g'(z_1)|}{|g'(z_2)|} \leq 8 \frac{z_2}{z_1} \frac{\ln( z_1/R)}{\ln( z_2/R)}\frac{|g(z_1)|^2}{|g(z_2)|^2}.
\end{align}

Now we proceed in a similar way as in the proof of \cite[Lemma~2.6]{S}.
If $|g(z_1)| > 1$, then \eqref{eq:2} implies
\begin{multline*}
\frac{|g^*(z_1)|}{|g^*(z_2)|} \geq \frac{1}{2} \frac{1+z_1^2}{z_1} \frac{z_2}{1+z_2^2} \left(\frac{\ln( z_2/R)}{\ln( z_1/R)}\right)^3
\frac{|g(z_1)|^2}{|g(z_1)|^2+1} \frac{|g(z_2)|^2+1}{|g(z_2)|^2}\\
> \frac{1}{8}\frac{z_1}{z_2} \left(\frac{\ln( z_2/R)}{\ln( z_1/R)}\right)^3 \geq c_4 \frac{z_1}{z_2} \left(\frac{\ln z_1}{\ln z_2}\right)^{-3},
\end{multline*}
and if $|g(z_1)| \leq 1$, then \eqref{eq:1} gives
\begin{multline*}
\frac{|g^*(z_1)|}{|g^*(z_2)|} \geq \frac{1}{2} \frac{1+z_1^2}{z_1} \frac{z_2}{1+z_2^2} \left(\frac{\ln( z_2/R)}{\ln( z_1/R)}\right)^3
\frac{|g(z_2)|^2+1}{|g(z_1)|^2+1}\\
> \frac{1}{8}\frac{z_1}{z_2} \left(\frac{\ln( z_2/R)}{\ln( z_1/R)}\right)^3 \geq c_4 \frac{z_1}{z_2} \left(\frac{\ln z_1}{\ln z_2}\right)^{-3},
\end{multline*}
where the constant $c_4 > 0$ depends only on $R,L$. 
Similarly, if $|g(z_2)| > 1$, then by \eqref{eq:2},
\[
\frac{|g^*(z_1)|}{|g^*(z_2)|} \leq 8 \frac{1+z_1^2}{z_1} \frac{z_2}{1+z_2^2} \frac{\ln( z_1/R)}{\ln( z_2/R)}
\frac{|g(z_1)|^2}{|g(z_1)|^2+1} \frac{|g(z_2)|^2+1}{|g(z_2)|^2}
< 32\frac{z_1}{z_2} \frac{\ln( z_1/R)}{\ln( z_2/R)} \leq c_5 \frac{z_1}{z_2} \frac{\ln z_1}{\ln z_2},
\]
and if $|g(z_2)| \leq 1$, then by \eqref{eq:1},
\[
\frac{|g^*(z_1)|}{|g^*(z_2)|} \leq 8 \frac{1+z_1^2}{z_1} \frac{z_2}{1+z_2^2} \frac{\ln( z_1/R)}{\ln( z_2/R)}
\frac{|g(z_2)|^2+1}{|g(z_1)|^2+1}
< 32\frac{z_1}{z_2} \frac{\ln( z_1/R)}{\ln( z_2/R)} \leq c_5 \frac{z_1}{z_2} \frac{\ln z_1}{\ln z_2},
\]
where the constant $c_5 > 0$ depends only on $R,L$. 

The second assertion of the lemma is obvious from the construction (e.g.~in the case $a = \infty$ the branches $g_j$, $j = 1,2$ can be defined on domains of the form $V \setminus \{z: \Arg(z) = t_j\}$ for suitable $t_1 \neq t_2$). 
\end{proof}

\begin{cor}\label{cor:log}
Let $U$ be a logarithmic tract of $f: U \to V$ over $a\in\clC$, where $V = \DD(a, r) \setminus\{a\}$ for some $0 < r < 1$. Then for every $z \in V$ there exists $w \in f^{-1}(z)$, such that $|w| \to +\infty$ and $|f^*(w)| \to 0$ as $z \to a$.
\end{cor}
\begin{proof} By diminishing $V$ if necessary, we can assume that $0 \notin U$. Fix a point $z_0 \in V$ with $d(z_0, a) < r/2$ and a branch $g$ of $f^{-1}$ near $z_0$. Take an arbitrary $z \in V$ with $d(z, a) \leq d(z_0, a)$. Then by the Spherical Distortion Theorem for logarithmic tracts for $z_1 = z$, $z_2 = z_0$, there exists $c_1 = c_1(z_0)$ and a branch $g$ of $f^{-1}$ in a neighbourhood of $z$ such that 
\[
|g^*(z)| \geq \frac{c_1}{d(z,a) (\ln(1/d(z,a))^3} \xrightarrow[z\to a]{}+\infty.
\]
This implies that setting $w = g(z)$, we have $|f^*(w)| \to 0$ as $z \to a$.
By the definition of a logarithmic tract, $|w| \to +\infty$ as $z \to a$. This proves the corollary.
\end{proof}

\begin{rem}\label{rem:cover} Looking at the proofs of the above results it is easy to see that in fact they are valid not only for maps $f$ defining a logarithmic tract over $a \in \clC$, but also for any holomorphic universal covering $f$ from a simply connected domain $U \subset \C$ onto $V = \DD(a, r) \setminus \{a\}$. 
\end{rem}

\begin{lem}\label{lem:diam}
Let $f$ be a transcendental meromorphic map. Then there exists a
sequence $d_n \to 0$, such that for every $n > 0$, every $z\in J(f)
\setminus \PP_n(f)$, every component of 
\[
f^{-n}(\DD(z, d(z, \PP_n(f) \cup \{\infty\})/2))
\]
has spherical diameter smaller than $d_n$. 
\end{lem}

\begin{proof}
Note that all branches of $f^{-n}$ are defined on $\DD(z, d(z, \PP_n(f) \cup \{\infty\}))$. 
Suppose the assertion does not hold. Then there exist $\delta > 0$, a
subsequence $n_k$, points $z_k \in J(f)$ and branches $g_k$ of
$f^{-{n_k}}$ on 
$\DD(z_k, d(z_k, \PP_{n_k}(f) \cup \{\infty\}))$, such that
$\diam_{sph} g_k(D_k) \geq \delta$ for every $k$, where
\[
D_k = \DD(z_k, d(z_k, \PP_{n_k}(f) \cup \{\infty\})/2).
\]
We have $\diam_{sph} D_k \leq\diam_{sph}\clC/2$ and $\diam_{sph}(\clC
\setminus g_k(D_k)) \geq \diam_{sph}(\Sing(f) \cup \{\infty\}) > 0$.
Hence, by the Spherical Koebe Distortion Theorem, $g_k(D_k)$ contains
a spherical disc of radius $\delta' > 0$ (independent of $k$).
Passing to a subsequence, we can assume that $z_k \to z$ for some $z
\in \clC$ and there is a spherical disc $D \subset \bigcap_k
g_k(D_k)$. Then 
\[
f^{n_k}(D) \subset \DD(z, 3d(z, \PP_{n_k}(f) \cup \{\infty\})/4)
\subset \DD(z, 3d(z, \Sing(f) \cup \{\infty\})/4)
\]
for large $k$, so the family $\{f^{n_k}|_D\}_{k>0}$ is normal. On the
other hand, since $z_k \in J(f)$, using the invariance of the Julia
set and the spherical Koebe Distortion Theorem we see that $D$
contains a point from $J(f)$, which gives a contradiction.
\end{proof}

Recall that a (Picard) exceptional value of $f$ is a point $a \in \clC$ such that $\bigcup_{n=0}^\infty f^{-n}(a)$ is finite. By Picard's Theorem, a meromorphic map has at most two exceptional values. Unlike the case of rational maps, for transcendental meromorphic maps the exceptional values can be contained in the Julia set. An exceptional value $a$ is called omitted, if $f^{-1}(a) = \emptyset$.

\begin{rem}\label{rem:excep}By Iversen's Theorem (see e.g.~\cite{CL}), every exceptional value $a$ is an asymptotic value, so $f$ has a singularity over $a$. If $a$ is an isolated point of $\Sing(f)$ (e.g.~if $f \in \SSS$), then all singularities over $a$ are logarithmic. If $f$ is entire, then it has at most one finite exceptional value, which is necessarily omitted.
\end{rem}

\begin{lem} \label{lem:compact} Let $f$ be a transcendental meromorphic map. Then for every compact set $K \subset J(f)$ such that $K$ does not contain exceptional values of $f$ and for every open set $U \subset \C$ intersecting $J(f)$ there exists $m \geq 0$ and $c > 0$ such that for every $z \in K$ 
we have 
\[
|(f^m)^*(w)| < c
\]
for some $w \in f^{-m}(z) \cap U$.
\end{lem}
\begin{proof} Since $U$ intersects $J(f)$, we have $f^m(J(f) \cap U) \supset K$ for a large $m$ (see e.g.~\cite{Ber}). Take a point $z \in K$. Then there exists a point $w \in J(f) \cap U$, such that $f^m(w) = z$. Hence, for some small neighbourhood $U_z$ of $w$, the map $f^m$ is defined on $U_z$, $f^m(U_z)$ is a neighbourhood of $z$ and $|(f^m)^*| < c_z$ on $U_z$ for some constant $c_z$. Since $\{f^m(U_z)\}_{z\in K}$ is an open cover of $K$, choosing a finite subcover we get the assertion.
\end{proof}

\section{Pressure for meromorphic maps and GPS points}\label{sec:pressure}

Let $f:\C\to\clC$ be a transcendental meromorphic map. For $z \in \clC$ and $t > 0$ denote by $S_n(t,z)$ the sum
\[
S_n(t,z)=\sum_{w\in f^{-n}(z)}|(f^n)^*(w)|^{-t}.
\]
For a set $A \subset \clC$ we will write
\[
S_n^A(t, z) = \sum_{w \in f^{-n}(z) \cap A} |(f^n)^*(w)|^{-t}.
\]

\begin{defn} Let $z \in \clC$ and $t > 0$. We define the lower and upper topological pressure for $f$ at the point $z$ as
\[
\underline{P} (f, t, z)=\liminf_{n\to\infty}\frac 1 n\ln S_n(t,z), \qquad \overline{P} (f, t, z)=\limsup_{n\to\infty}\frac 1 n\ln S_n(t,z).
\]
Note that the values $\pm\infty$ are not excluded. In the formulations of the results concerning pressure, we will usually include both the finite and the infinite case, considering the standard order and topology in $\overline{\R} = \R\cup\{\pm\infty\}$. 
\end{defn}

If the lower and upper pressures coincide, i.e.~if there exists the limit
\[
P(f, t, z)=\lim_{n\to\infty}\frac 1 n\ln S_n(t,z),
\]
then we call $P(f, t, z)$ the topological pressure for $f$ at the point $z$.
Denoting the pressure, we will often omit the symbol $f$. 

\begin{defn}
Let $f$ be a meromorphic map. We will say that a point $z \in \C$ is a GPS point (Good Pressure Starting point) if 
\[
\ln d(z, \PP_n(f)) = o(n),
\]
i.e. $d(z, \PP_n(f)) > e^{-a_n}$ for large $n$, where $a_n > 0$, $\frac{a_n}{n} \to 0$ (the sequence $a_n$ may depend on $z$).
\end{defn}

The following proposition is straightforward.

\begin{prop}\label{prop:GPS}

\

\begin{itemize}
\item 
Every preimage of a GPS point is GPS.
\item
Every GPS point is in $\C\setminus \PP(f)$.
\item
Every point $z \in \C \setminus \overline{\PP(f)}$ is GPS.
\qed
\end{itemize}
\end{prop}

Recall that a conformal repeller $X$ is transitive, if for all non-empty sets $U, V$ open in $X$ we have $f^n(U) \cap V \neq \emptyset$ for some $n \geq 0$; $X$ is isolated, if there exists a neighbourhood $W$ of $X$, such that for every $z \in W\setminus X$ there exists $n >0$ with $f^n(z) \notin W$.

\begin{defn}
Let $P_{hyp}(t)$ be the supremum of the pressures $P(f|_X, t)$ over all transitive isolated conformal repellers $X \subset J(f)$.
\end{defn}

\section{Pressure for maps in $\SSS$}\label{sec:pressure-S}

\begin{prop}\label{prop:E}
If $f \in \SSS$, then the set of non-GPS points has Hausdorff dimension $0$. 
\end{prop}
\begin{proof}
Suppose $f \in \SSS$ and let
\[
E = \bigcap_{k=0}^\infty \bigcup_{n=k}^\infty \DD(f^n(\Sing(f)), e^{-\sqrt{n}}) \cup \PP(f) \cup \{\infty\},
\]
Since $f \in \SSS$, for every $\delta > 0$ and $k >0$ the set $E$ can be
covered by a countable number of spherical discs $D_j$ such that
$\sum_j (\diam_{sph} D_j)^\delta < c \sum_{n=k}^\infty n
e^{-\delta\sqrt{n}} \to 0$ as $k \to \infty$. This shows that
$\dim_H E = 0$.

Take $z \in \C \setminus E$. Then there exists $k$ such that $d(z, f^n(\Sing(f))) \geq  e^{-\sqrt{n}}$ for every $n \geq k$. Moreover, $d(z, f^n(\Sing(f))) > 0$ for $0 \leq n < k$. This implies that 
$d(z, \PP_n(f))) \geq e^{-\sqrt{n}}$ for large $n$, which shows that $z$ is a GPS point. Hence, the set of non-GPS points is contained in $E$. 
\end{proof}

\begin{thm}\label{thm:pressure}
Let $f$ be a transcendental meromorphic function in the class $\SSS$. Then $\underline{P}(t,z)$ and $\overline{P}(t,z)$ do not depend on $z$ within the set of all GPS points.
\end{thm}

By the above theorem, for maps in class $\SSS$ we can define respectively by $\underline{P}(t)$ and $\overline{P}(t)$ the common values 
$\underline{P}(t,z)$ and $\overline{P}(t,z)$ for all GPS points $z$.

The proof of Theorem~\ref{thm:pressure} follows \cite{P}. In particular, we use the following lemma from this paper.

\begin{lem}[\mbox{\cite[Lemma 3.1]{P}}]\label{lem:points}
There exists $C > 0$ such that for every finite set $W$ of points in $\clC$ and $0<r<1/2$, for every $z_1, z_2 \in \clC \setminus \DD(W, r)$ there exists a sequence of spherical discs $D_1 = \DD(q_1, \rho_1), \ldots, D_k = \DD(q_1, \rho_k)$, such that $\bigcup_{j=1}^k D_j$ is connected, $z_1\in D_1$, $z_2 \in D_k$, for every $j = 1, \ldots, k$ the spherical disc $\DD(q_j, 2\rho_j)$ is disjoint from $W$ and 
\[
\begin{array}{l r}
k \leq C\sqrt{\# W}\sqrt{\ln 1/r} & \text{if } \# W \geq \ln 1/r,\\
k \leq C\ln 1/r & \text{if } \# W < \ln 1/r.  
\end{array}
\]
\end{lem}

Using the above lemma, we show the following.

\begin{lem} \label{lem:S_n} If $f \in \SSS$, then for all GPS points $z_1, z_2$,
\[
S_n(t, z_1) \leq e^{o(n)} S_n(t, z_2).
\]
\end{lem}
\begin{proof}
By definition, there exists a sequence $0 < a_n = o(n)$, such that 
$d(z_j, \PP_n(f)) > e^{-a_n}$ for large $n$ and $j = 1,2$. 
Hence, for large $n$ the points $z_1, z_2$ satisfy the assumptions of Lemma~\ref{lem:points} for $W = \PP_{n+1}(f)$, $r = e^{-a_{n+1}}$. Note that $\# W < c_1(n+1)$ for some $c_1$ and $\ln 1/r = a_{n+1}$. Hence, Lemma~\ref{lem:points} implies that there exist $k$ suitable discs $D_j = \DD(q_j, \rho_j)$, where 
\[
k \leq C \max(\sqrt{(n+1)a_{n+1}}, a_{n+1}) = O(a_{n+1}) = o(n),
\]
such that all branches of $f^{-(n+1)}$ are defined on $\DD(q_j, 2\rho_j)$. Moreover, for every branch $g$ of $f^{-n}$ on $D_j$, we have
$\diam_{sph} D_j < \diam_{sph} \clC / 2$ and $\clC \setminus g(D_j)$ contains $\Sing(f) \cup \{\infty\}$. Hence, by the Spherical Koebe Distortion Theorem, 
\begin{equation}\label{eq:g^*}
\frac{|g^*(z)|}{|g^*(z')|} < c_2
\end{equation}
for all $z, z' \in D_j$, where $c_2$ does not depend on $n, j, g$. This implies
\[
S_n(t, z) \leq c_2^t S_n(t, z').
\]
Since $\bigcup_{j=1}^k D_j$ is connected and $z_1\in D_1$, $z_2 \in D_k$, we have
\[
S_n(t, z_1) \leq c_2^{kt} S_n(t, z_2) = e^{o(n)} S_n(t, z_2).
\]
\end{proof}

\begin{proof}[Proof of Theorem~\rm\ref{thm:pressure}]
Take two GPS points $z_1, z_2$. Then by Lemma~\ref{lem:S_n}, 
\[
\frac{1}{n} \ln S_n(t, z_1) \leq \frac{1}{n} \ln S_n(t, z_2) + \frac{o(n)}{n},
\]
which implies $\underline{P}(t, z_1) \leq \underline{P}(t, z_2)$ and $\overline{P}(t, z_1) \leq \overline{P}(t, z_2)$. By symmetry, $\underline{P}(t, z_1) = \underline{P}(t, z_2)$ and $\overline{P}(t, z_1) = \overline{P}(t, z_2)$.
\end{proof}

The following theorem is the main result of this section. 

\begin{thm}\label{thm:pressure-S}
Let $f$ be a transcendental meromorphic map in the class $\SSS$. Then 
\[
\underline{P}(t) = \overline{P}(t) = P_{hyp}(t) \quad \text{for every }  t  > 0.
\]
\end{thm}

Before the proof of the theorem, we state two propositions.

In \cite{S3} (see also \cite{Ber2}) it was proved that $\dim_H(J(f)) > 0$ for all transcendental meromorphic maps $f$. Since the proof is done by constructing an invariant hyperbolic Cantor set of positive Hausdorff dimension, in fact we get:

\begin{prop}\label{prop:hyp} If $f:\C\to\clC$ is meromorphic, then $\dim_{hyp}(J(f)) > 0$. In fact, for every open set $U$ intersecting $J(f)$
there exists a transitive isolated Cantor conformal repeller $X \subset U \cap J(f)$ of positive Hausdorff dimension. 
\end{prop}

Recall that for a set $A \subset \clC$ we write
\[
S_n^A(t, z) = \sum_{w \in f^{-n}(z) \cap A} |(f^n)^*(w)|^{-t}.
\]

The proof of the following proposition repeats the proof of
\cite[Proposition~2.1]{PRS}.

\begin{prop}\label{prop:S^K} Let $f \in \SSS$. Then there exists a GPS point $z_0 \in J(f)$, such that for every $t > 0$,
\[
P_{hyp}(t) \geq \sup_K \limsup_{n\to\infty} \frac{1}{n} \ln S_n^K(t, z_0),
\]
where the supremum is taken over all compact sets $K \subset J(f)$, which do not contain exceptional values of $f$. 
\end{prop}
\begin{proof}
By Propositions~\ref{prop:E} and~\ref{prop:hyp}, there exists an isolated transitive conformal repeller $X \subset J(f)$ and a GPS point $z_0 \in X$. 
Then for sufficiently large $n$, we have $d(z, \PP_n(f)) > e^{-a_n}$
for some sequence $a_n > 0$ with $a_n/n \to 0$, so all branches of $f^{-n}$ are defined on the spherical disc $\DD(z_0, e^{-a_n})$. Obviously, we can assume that $a_n \to \infty$ as $n \to \infty$.

By the definition of a repeller, there exist constants $\delta, c_1 > 0$, $Q_1,Q_2 > 1$, such that for every $z \in X$ and $l > 0$, there exists a branch $g^l$ of $f^{-l}$ defined on $\DD(f^l(z), \delta)$, with $g^l(f^l(z)) = z$ and
\begin{equation}\label{eq:f^l*}
Q_1^{-l}/c_1 < |(g^l)^*| < c_1 Q_2^{-l}. 
\end{equation}
Fix a large constant $M$ and take
\[
l_n = [M a_{2n}]
\]
for large $n$. Let 
\begin{equation}\label{eq:D}
D_n = \DD(f^{l_n}(z_0), \delta)
\end{equation}
and let $g^{l_n}$ be the branch of $f^{-l_n}$ on $D_n$ such that $g^{l_n}(f^{l_n}(z_0)) = z_0$. Then by \eqref{eq:f^l*}, we have $|(g^{l_n})^*| < c_1 Q_2^{-l_n}$ on $D_n$, so if $M$ is chosen sufficiently large, then 
\begin{equation}\label{eq:g^l}
g^{l_n}(D_n) \subset \DD(z_0, e^{-a_{2n}}/2).
\end{equation}
Take a compact set $K \subset J(f)$, which does not contain the exceptional values of $f$. By Lemma~\ref{lem:compact} and the compactness of $X$, there exists a bounded sequence $m_n$ such that for every $z \in K$ we can find a point $w_z \in \DD(f^{l_n}(z_0), \delta/2)$ with $f^{m_n}(w_z) = z$ and
\begin{equation}\label{eq:f^m*}
|(f^{m_n})^*(w_z)| < c_2
\end{equation}
for some $c_2 > 0$. Since $z_0$ is a GPS point and the sequence $m_n$ is bounded, for sufficiently large $n$ all the inverse branches of $f^{-(n+m_n)}$ are defined in $\DD(z_0, e^{-a_{2n}})$. Hence, for every $z \in f^{-n}(z_0) \cap K$ there exists a branch $g_z$ of 
$f^{-n}$ on $\DD(z_0, e^{-a_{2n}})$, such that $g_z(z_0) = z$ and a branch $\tilde g_{w_z}$ of $f^{-m_n}$ on $g_z(\DD(z_0, e^{-a_{2n}}))$, such that $\tilde g_{w_z}(z) = w_z$ for $w_z$ from \eqref{eq:f^m*}.

Let 
\[
N = N(n) = l_n+n+m_n
\]
and consider the family 
\[
\F_n = \{h_z\}_{z \in f^{-n}(z_0) \cap K}
\]
of inverse branches of $f^{-N}$ defined on $D_n$ from \eqref{eq:D}, where 
\[
h_z = \tilde g_{w_z}\circ g_z \circ g^{l_n}.
\]
By \eqref{eq:g^l} and Lemma~\ref{lem:diam}, the spherical diameter of every set $h_z(D_n)$ for $h_z \in \F_n$ is smaller than $\delta/4$, if $n$ is sufficiently large. Since $w_z \in \DD(f^{l_n}(z_0), \delta/2)$, this implies 
\begin{equation}\label{eq:h_z(D)}
\overline{h_z(D_n)} \subset \overline{\DD(f^{l_n}(z_0), 3\delta/4)} \subset D_n.
\end{equation}

Hence, the family $\F_n$ forms a conformal iterated function system on $D_n$ (it is finite since $K$ is compact) and its limit set $\Lambda_n$ is a transitive isolated $f^N$-invariant conformal Cantor repeller contained in $\bigcup_{z \in f^{-n}(z_0) \cap K}h_z(D_n)$.
We have
\[
\sum_{z \in f^{-n}(z_0):\;  h_z \in \F_n} |(f^n)^*(z)|^{-t} = S_n^K(t, z_0), 
\]
which together with \eqref{eq:f^l*} and \eqref{eq:f^m*} implies
\[
\sum_{h_z \in \F_n} |h_z^*(f^{l_n}(z_0))|^t \geq \frac{c_1^t}{c_2^t Q_1^{l_nt}} S_n^K(t, z_0).
\]
We have $\diam_{sph}(D_n) = r_1$ and, by \eqref{eq:h_z(D)}, $\diam_{sph}(\clC \setminus h_z(D_n)) > \diam_{sph}(\clC \setminus D_n) = r_2 > 0$ for $r_1, r_2$ independent of $n$, so we can use the Spherical Koebe Distortion Theorem to obtain
\begin{equation}\label{eq:f^n*}
\sum_{h_z \in \F_n} |h_z^*(v)|^t \geq \frac{c_3}{Q_1^{l_nt}} S_n^K(t, z_0) 
\end{equation}
for every $v \in \Lambda_n$ and a constant $c_3 > 0$ (depending on $t$). Let 
\[
Y_n = \bigcup_{j = 0}^{N-1} f^j(\Lambda_n).
\]
Then $Y_n$ is a transitive isolated $f$-invariant conformal repeller, so the pressure function of $f|_{Y_n}$ is equal to 
\begin{multline*}
P(f|_{Y_n}, t) = \lim_{k\to\infty} \frac{1}{k} \ln S^{Y_n}_k(t,v) \\=
\lim_{k\to\infty} \frac{1}{Nk} \ln S^{\Lambda_n}_{Nk}(t,v)
= \lim_{k\to\infty} \frac{1}{Nk} \ln \sum_{\substack{h_j \in \F_n\\ \text{for } j = 1, \ldots, k}} |(h_k\circ \cdots \circ h_1)^*(v)|^t
\end{multline*}
for every $v \in \Lambda_n$. Hence, applying \eqref{eq:f^n*} to $h_j\circ \cdots \circ h_1(v)$ for $j = 0,\ldots, k-1$ we get
\[
P(f|_{Y_n}, t) \geq \lim_{k\to\infty} \frac{1}{Nk} \ln \left(\frac{c_3}{Q_1^{l_nt}} S_n^K(t, z_0)\right)^k = \frac{1}{N} \ln S_n^K(t, z_0) + \frac{\ln c_3 - l_nt \ln Q_1}{N}.
\]
Since $l_n = O(a_{2n}) = o(n)$ and $N = n + O(a_{2n}) = n + o(n)$ as $n \to \infty$, this implies
\[
P_{hyp}(t) \geq \limsup_{n\to\infty} P(f|_{Y_n}, t) \geq \limsup_{n\to\infty} \frac{1}{n} \ln S_n^K(t, z_0).
\]
\end{proof}

\begin{proof}[Proof of Theorem~\rm\ref{thm:pressure-S}]

First we prove $P_{hyp}(t) \leq \underline{P}(t)$. Take a transitive isolated conformal repeller $X \subset J(f)$ and a point $z \in X$. Then
\begin{equation}\label{eq:f|X}
P(f|_X, t) = \lim_{n\to\infty} \frac{1}{n} \ln S_n^X(t, z)
\end{equation}
(see e.g.~\cite{PU}). By the properties of the repeller, there exists $\delta > 0$, such that for every $n >0$, all branches $g$ of $f^{-n}$ with $g(z) \in X$ are defined on $\DD(z, \delta)$. Proposition~\ref{prop:E} implies that there exists a GPS point $z_0 \in \DD(z, \delta/2)$. By the Spherical Koebe Distortion Theorem, the spherical distortion of $g$ on $\DD(z, \delta/2)$ is universally bounded. This together with \eqref{eq:f|X} gives
\[
P(f|_X, t) \leq \underline{P}(t,z_0),
\]
which shows $P_{hyp}(t) \leq \underline{P}(t)$.

Now we prove $P_{hyp}(t) \geq \overline{P}(t)$. The proof will be split into three cases:

\begin{description}
\item[\rm\emph{Case $1$}] $f$ has no poles.
\item[\rm\emph{Case $2$}] $f$ has a pole and all poles are exceptional values of $f$.
\item[\rm\emph{Case $3$}] $f$ has a pole, which is not an exceptional value of $f$.
\end{description}

\subsubsection*{Case $1$} 

Fix $t > 0$ and an arbitrary $\varepsilon > 0$. Let $Q= e^{\overline{P}(t)- \varepsilon}$ if $\overline{P}(t)$ is finite, or set $Q$ to be an arbitrary fixed number if $\overline{P}(t)$ is infinite. We will show that $P_{hyp}(t) \geq \ln Q$, which will end the proof.

Since $f\in\SSS$ and $f$ has no poles, there exists $V = \{z\in\C: |z| > R\}$ for some $R > 1$ such that $f^{-1}(V)$ is a union of logarithmic tracts of $f$ over $\infty$. Enlarging $R$, we can assume that all these tracts do not contain $0$. 

By Corollary~\ref{cor:sta}, we have $S_1(t,v) \to +\infty$ when $|v| \to +\infty$, so enlarging $R$ we can assume that for $|v| = 2R$ we have
\[
S_1(t,v) > \frac{3}{c^t} Q
\]
for $c$ from Corollary~\ref{cor:sta} (for $R$ as above and $L = 2$). Hence, there exists $M_0 > 0$ such that 
\begin{equation}\label{eq:S_1>}
S_1^{\clD(M_0)}(t,v) > \frac{2}{c^t} Q.
\end{equation}

Since $f$ is entire, it has at most one exceptional value $a \in J(f)$ and $a$ is omitted. If such a value exists, let $V_0$ be a small spherical disc centred at $a$. By Remark~\ref{rem:excep}, the set $f^{-1}(V_0)$ is a union of logarithmic tracts of $f$. 

Choosing $R$ large enough, we can assume
\begin{equation}\label{eq:V_0}
f(V_0) \subset \D(2R).
\end{equation}
Take a GPS point $z_0 \in J(f)$, for which Proposition~\ref{prop:S^K} holds and let 
\[
W = \overline{\D}(2R)\setminus V_0.
\]
Suppose that there exists $q > 0$ such that $S^W_n(t,z_0) > q Q^n$ for infinitely many $n$. Then, since $W$ is compact and does not contain the exceptional values of $f$, by Proposition~\ref{prop:S^K} for $K=W\cap J(f)$ we have $P_{hyp}(t) \geq \ln Q$, which is the assertion we want to show. Hence, we can assume that for every $n$, 
\begin{equation}\label{eq:S^K<}
S^W_n(t, z_0) < \varepsilon_n Q^n,
\end{equation}
where $\varepsilon_n \to 0$ as $n\to\infty$. 

By the definition of $Q$, we can take a subsequence $n_j \to \infty$, such that
\[
S_{n_j}(t, z_0) \geq Q^{n_j}.
\]
By \eqref{eq:S^K<}, this implies $S^{\C\setminus\overline{\D}(2R)}_{n_j}(t, z_0) \geq Q^{n_j}/3$ or $S^{V_0}_{n_j}(t, z_0) \geq Q^{n_j}/3$ for large $j$. In the second case, by Corollary~\ref{cor:log}, we have $S^{\C\setminus\overline{\D}(2R)}_{n_j + 1}(t, z_0) \geq Q^{n_j + 1}/3$ (if $V_0$  is sufficiently small). We conclude that in both cases, for arbitrarily large numbers $m$ (equal to $n_j$ or $n_j+1$) we have 
\[
S^{\C\setminus\D(2R)}_m(t, z_0) \geq  \frac{Q^m}{3}.
\]
Fix such a number $m$ and take $M > 2R$ so large that 
\begin{equation}\label{eq:M}
M > \tilde c M_0 \left(\frac{\ln M}{\ln(2R)}\right)^{4\pi},
\end{equation}
where $\tilde c$ is the constant from Corollary~\ref{cor:|g|_infty} (for $R$ as above and $L=2$) and
\begin{equation}\label{eq:S_m>}
S^A_m(t, z_0) \geq\frac{Q^m}{6},
\end{equation}
where
\[
A = \overline{\D}(M)\setminus\D(2R).
\]

Take an arbitrary $z\in A$. Then, using successively Corollary~\ref{cor:|g|_infty} (for $z_1=z$, $z_2=v$), \eqref{eq:M}, Corollary~\ref{cor:sta} and \eqref{eq:S_1>}, we obtain
\[
S^{\overline{\D}(M)}_1(t, z) \geq \sum_{\substack{g \in f^{-1},\\|g(v)| \leq M_0}} |g^*(z)|^t \geq c^t \frac{(\ln(2R))^{3t}}{(2R)^t} \frac{|z|^t}{(\ln|z|)^{3t}} S_1^{\clD(M_0)}(t,v)
\geq 2Q \frac{(\ln(2R))^{3t}}{(2R)^t} \frac{|z|^t}{(\ln|z|)^{3t}}.
\]
Since the function $x \mapsto x^t/(\ln x)^{3t}$ is increasing for $x > 2R$ (provided $R$ was chosen large enough), this yields
\begin{equation}\label{eq:S_1>2Q}
S^{\overline{\D}(M)}_1(t, z) \geq 2Q \qquad \text{for every } z \in A.
\end{equation}
By \eqref{eq:S_m>} and \eqref{eq:S_1>2Q}, 
\[
\sum_{\substack{w\in f^{-(m+1)}(z_0)\cap \clD(M),\\f(w)\in A}}
|(f^{m+1})^*(w)|^{-t} \geq \frac{Q^{m+1}}{3},
\]
which together with \eqref{eq:V_0} and \eqref{eq:S^K<} implies 
\begin{align*}
S_{m+1}^A(t, z_0) &\geq \sum_{\substack{w\in f^{-(m+1)}(z_0)\cap A,\\f(w)\in A}}
|(f^{m+1})^*(w)|^{-t}\\ &= \sum_{\substack{w\in f^{-(m+1)}(z_0)\cap \clD(M),\\f(w)\in A}}|(f^{m+1})^*(w)|^{-t} - 
\sum_{\substack{w\in f^{-(m+1)}(z_0)\cap W,\\f(w)\in A}}
|(f^{m+1})^*(w)|^{-t}\\ &\geq \left(\frac{1}{3}-\varepsilon_{m+1}\right)Q^{m+1} > \frac{Q^{m+1}}{6}
\end{align*}
(if $m$ was chosen large enough). The latter inequality is the same as \eqref{eq:S_m>}, with $m$ replaced by $m+1$. Hence, by induction, we get
\[
S_n^A(t, z_0) \geq \frac{Q^n}{6}
\]
for every $n \geq m$. Using Proposition~\ref{prop:S^K} for the set $K' = A \cap J(f)$, we obtain $P_{hyp}(t) \geq \ln Q$, which ends the proof. 

\subsubsection*{Case $2$} Suppose that $f$ has a pole and all poles are exceptional values of $f$. Then $\infty$ is also an exceptional value, so in fact there is only one pole $a$ and $a$ is omitted. This implies that $f$ is a self-map of the punctured plane $\C\setminus\{a\}$. By a change of coordinates, we can assume $a = 0$. 

Let $V = \{z\in\C: |z| > R\}$ and $V' = \D(1/R) \setminus \{0\}$ for a large $R > 0$. By Remark~\ref{rem:log}, the set $f^{-1}(V)$ is a union of logarithmic tracts of $f$ over $\infty$ (contained in a small neighbourhood of $\infty$) and a small punctured neighbourhood $W$ of $0$, such that $f$ on $W$ is a finite degree covering. Similarly, by Remark~\ref{rem:excep}, the set $f^{-1}(V')$ is a union of logarithmic tracts of $f$ over $0$ (contained in a small neighbourhood of $\infty$). This implies that for $R$ large enough, the sets $f^{-2}(V)$, $f^{-2}(V')$ are unions of logarithmic tracts of $f^2$ over $0$ or $\infty$ and some simply connected domains contained in a small punctured neighbourhood of $0$, such that $f^2$ is a universal covering on each of these domains.

Now we virtually repeat the proof in Case~1, replacing $S_1(t,z)$ by $S_2(t, z)$. Note that by Remark~\ref{rem:cover}, we can use the distortion results from Section~\ref{sec:lem} for all components of $f^{-2}(V)$, $f^{-2}(V')$ (we shall not repeat this remark in the sequel).

Define $Q$ as in Case~1. By Corollary~\ref{cor:sta}, enlarging $R$ we can assume that for some $v$ with $|v| = 2R$ we have
\[
S_2(t, v), S_2(t, 1/v) > \frac{3}{c^t} Q^2,  
\]
where the Spherical Distortion Theorem for logarithmic tracts and Corollary~\ref{cor:sta} hold with the constant $c$ (for $R$ as above, $L = 2$,  $r = 1/R$, $\lambda = 1/2$). Hence, we can take $M_0 > 0$ such that 
\begin{equation}\label{eq:S_1>'}
S_2^{\clD(M_0)\setminus \D(1/M_0)}(t,v), 
S_2^{\clD(M_0)\setminus \D(1/M_0)}(t,1/v)
> \frac{2}{c^t} Q^2.
\end{equation}
As in Case~1, we take a GPS point $z_0 \in J(f)$ for which Proposition~\ref{prop:S^K} holds. Let 
\[
W = \overline{\D}(2R)\setminus \D(1/(2R)).
\]
By Proposition~\ref{prop:S^K}, we can assume that
\begin{equation}\label{eq:S^K<'}
S^W_n(t, z_0) < \varepsilon_n Q^n,
\end{equation}
where $\varepsilon_n \to 0$ as $n\to\infty$. 

Taking a subsequence $n_j \to \infty$, such that
\[
S_{n_j}(t, z_0) \geq Q^{n_j},
\]
and using \eqref{eq:S^K<'}, we have $S^{\C\setminus W}_{n_j}(t, z_0) \geq Q^{n_j}/2$ for large $j$. Set $m = n_j$ and take $M > 2R$ so large that 
\begin{equation}\label{eq:M'}
M > \tilde c M_0 \left(\frac{\ln M}{\ln(2R)}\right)^{4\pi},
\end{equation}
where $\tilde c$ is the constant from Corollary~\ref{cor:|g|_infty} (for $R$ as above and $L=2$) and 
\begin{equation}\label{eq:S_m>'}
S^A_m(t, z_0) \geq \frac{Q^m}{4},
\end{equation}
where
\[
A = (\overline{\D}(M)\setminus\D(2R))\cup (\overline{\D}(1/(2R))\setminus\D(1/M)).
\]

Take an arbitrary $z\in A$. As in Case~1, using successively Corollaries~\ref{cor:|g|_infty}-\ref{cor:|g|}, \eqref{eq:M'}, the Spherical Distortion Theorem for logarithmic tracts, Corollary~\ref{cor:sta} and \eqref{eq:S_1>'}, we obtain
\begin{equation}
\label{eq:S_1>2Q'}
S^{\clD(M) \setminus \D(1/M)}_2(t, z) \geq 2Q^2. 
\end{equation}
By \eqref{eq:S_m>'} and \eqref{eq:S_1>2Q'}, 
\[
S^{\clD(M) \setminus \D(1/M)}_{m+2}(t, z_0) \geq \frac{Q^{m+2}}{2},
\]
which together with \eqref{eq:S^K<'} implies
\[
S_{m+2}^A(t, z_0) \geq \frac{Q^{m+2}}{4}.
\]
This is the same as \eqref{eq:S_m>'}, with $m$ replaced by $m+2$. By induction, we get
\[
S_{m+2n}^A(t, z_0) \geq \frac{Q^{m+2n}}{4}
\]
for every $n \geq 0$. Using Proposition~\ref{prop:S^K} for the set $K' = A \cap J(f)$, we obtain $P_{hyp}(t) \geq \ln Q$, which ends the proof.

\subsubsection*{Case $3$} 
Suppose now that $f$ has a pole $p$, which is not an exceptional value.
Let $W$ be a small neighbourhood of $p$ and let $V_0$ be the union of small 
spherical discs centred at the exceptional values $a_1, a_2$ of $f$ with $a_1, a_2 \in J(f)$ (if such values exist), such that $W \cap V_0 = \emptyset$. Let $V_1 = \{z\in\C: |z| > R\}$ for $R$ so large, that $W \cup V_0 \subset \C\setminus V_1$ and $f(W) \supset V_1$. By Remark~\ref{rem:excep}, $f$ has logarithmic singularities over $a_1$ and $a_2$. Let $U$ be the union of some logarithmic tracts of $f$ over $a_1$ and $a_2$, such that $f(U) = V_0\setminus \{a_1, a_2\}$. Diminishing $V_0$, we can assume $U \subset V_1$. 

Take a GPS point $z_0 \in J(f)$ for which Proposition~\ref{prop:S^K} holds and let
\[
K = J(f) \setminus (V_0 \cup V_1).
\]
We will show that for every $n$,
\begin{equation}\label{eq:S^K}
S^K_{n+i_n}(t,z_0) \geq q S_n(t,z_0)
\end{equation}
for some constant $q > 0$ and some $i_n \in \{0,1,2\}$. To see this, note that
for each $n$ we have $S_n^K(t, z_0) > S_n(t, z_0)/3$ or $S_n^{V_0}(t, z_0) \geq S_n(t, z_0)/3$ or $S_n^{V_1}(t, z_0) \geq S_n(t, z_0)/3$. If the first possibility holds, then \eqref{eq:S^K} is satisfied with $i_n = 0$. If the third possibility takes place, then (since the spherical derivative of $f$ is bounded on $W$), we have 
\begin{equation}\label{eq:third}
S_{n+1}^K(t, z_0) \geq S_{n+1}^W(t, z_0) \geq c_1 S_n^{V_1}(t, z_0) \geq c_2 S_n(t, z_0)
\end{equation}
for some constants $c_1, c_2 > 0$, so \eqref{eq:S^K} holds with $i_n = 1$. Finally, if the second possibility is satisfied, then by Corollary~\ref{cor:log}, 
\[
S_{n+1}^{V_1}(t, z_0) \geq S_{n+1}^U(t, z_0) \geq c_3 S_n^{V_0}(t, z_0) \geq c_4 S_n(t, z_0)
\]
for some constants $c_3, c_4 > 0$, so repeating the argument used in \eqref{eq:third}, we obtain
\[
S_{n+2}^K(t, z_0) \geq c_5 S_n(t, z_0)
\]
for some constant $c_5 > 0$, which gives \eqref{eq:S^K} with $i_n = 2$. In this way we have shown \eqref{eq:S^K}. Now by \eqref{eq:S^K} and Proposition~\ref{prop:S^K} we get $P_{hyp}(t) \geq \overline{P}(t)$. This ends the proof of the theorem.
\end{proof}

\section{Bowen's formula for maps in $\SSS$}\label{sec:bowen-S}

Theorem~\ref{thm:pressure-S} enables us to make the following definition:

\begin{defn}
If $f\in\SSS$, then the pressure of $f$ is defined as
\[
P(t) = P(f, t, z) = \lim_{n\to\infty} \frac{1}{n} \ln S_n(t, z)
\]
for any GPS point $z \in \C$.
\end{defn}

\begin{prop}\label{prop:pressure}
If $f \in \SSS$, then:
\begin{itemize}
\item 
$P(t) > -\infty$ for every $t > 0$.
\item 
the function $t \mapsto P(t)$ is non-increasing for $t \in (0, +\infty)$,
\item
the function $t \mapsto P(t)$ is convex $($and hence continuous$)$ for $t \in (t_0, +\infty)$, where $t_0 = \inf\{t > 0: P(t) \text{ is finite}\}$,
\item  $P(2) \leq 0$.
\end{itemize}
\end{prop}
\begin{proof} The first assertion follows directly from Theorem~\ref{thm:pressure-S} (the pressure of $f$ on a conformal repeller is finite). To prove the second one, take a GPS point $z$ and $0 < t_1 < t_2$. Let
\[
D_n = \DD(z, d(z, \PP_n(f)\cup\{\infty\})/2).
\]
By the definition of GPS points, $\diam_{sph} D_n =1/ e^{o(n)}$, and by the Spherical Koebe Distortion Theorem, the spherical distortion of every branch of $f^{-n}$ on $D_n$ is universally bounded. Hence, 
\begin{multline*}
S_n(t_2,z) \leq S_n(t_1, z) \sup_{w \in f^{-n}(z)} \frac{1}{|(f^n)^*(w)|^{t_2-t_1}} \\\leq  c_1 S_n(t_1, z) \left(\frac{\sup\{\diam_{sph} U: U \text{ is a component of } f^{-n}(D_n)\}}{\diam_{sph} D_n}\right)^{t_2-t_1}\\
\leq c_1 (e^{o(n)})^{t_2-t_1} S_n(t_1, z),
\end{multline*}
where $c_1 >0$ is independent of $n$. This implies that $P(t_2) \leq P(t_1)$. 

The third assertion follows by the H\"older inequality. 

To prove the fourth one, note that for every branch $g$ of $f^{-n}$ on $D_n$, the spherical area (denoted by $\Area$) of $g(D)$ satisfies
\[
\Area (g(D_n)) > c_2  \frac{|g^*(z)|^2}{e^{o(n)}}
\]
for a constant $c_2 > 0$. Since the sets $g(D_n)$ are disjoint for different branches $g$, this implies
\[
P(2) = \lim_{n\to\infty} \frac{1}{n}\ln \sum_g |g^*(z)|^2 \leq 
\lim_{n\to\infty} \frac{1}{n}\ln \left(e^{o(n)} \Area(\clC)\right) = 0.
\]
\end{proof}
 
\begin{defn} For $f \in \SSS$ let
\[
\delta(f) = \inf \{t > 0: P(t)  \leq 0\}.
\]
By Proposition~\ref{prop:pressure}, it is well-defined and $\delta(f) \leq 2$.
\end{defn}

As a corollary of Theorem~\ref{thm:pressure-S}, we immediately obtain the following version of Bowen's formula. 

\begin{bowen-S}
If $f$ is a transcendental meromorphic map in the class $\SSS$, then 
\[
\dim_H J_r(f)  = \dim_{hyp} J(f) = \delta(f).
\]
\end{bowen-S}
\begin{proof}
The equality between the Hausdorff dimension of $J_r(f)$ and hyperbolic dimension of $J(f)$ follows from \cite{R}. 

To show $\dim_{hyp} J(f) = \delta(f)$ using Theorem~\ref{thm:pressure-S}, it is enough to notice that the Hausdorff dimension of a transitive isolated conformal repeller $X \subset J(f)$ is, by the classical Bowen's formula, equal to the unique zero of the pressure 
\[
P(f|_X, t) = \lim_{n\to\infty} \frac{1}{n} \ln \sum_{w \in f^{-n}(z) \cap X} |(f^n)^*(w)|^{-t}
\]
for $z \in X$, which is a strictly decreasing function of $t$ (see e.g.~\cite{PU}). 
\end{proof}

\section{Pressure and Bowen's formula for maps in $\B$}\label{sec:pressure-B}

In this section we consider maps $f \in \B$, such that $J(f) \setminus \overline{\PP(f)} \neq \emptyset$. The results we prove are similar to the ones for the class $\SSS$, with analogous proofs. Hence, we only sketch them, indicating the differences compared to the case $f \in \SSS$. 

By Proposition~\ref{prop:GPS} and the fact that the preimages of points in $J(f)$ are dense in $J(f)$, we immediately obtain the following proposition.

\begin{prop}\label{prop:E'}
If $f$ is meromorphic and $J(f) \setminus \overline{\PP(f)} \neq \emptyset$, then $J(f) \setminus \overline{\PP(f)}$ is an open and dense backward-invariant subset of $J(f)$, contained in the set of GPS points. 
\qed
\end{prop}

\begin{defn} Let $f\in \B$. We will call $f$ exceptional, if there exists an exceptional value $a$ of $f$, such that $a \in J(f)$ and $f$ has a non-logarithmic singularity over $a$. Otherwise, we will say that $f$ is non-exceptional.
\end{defn}

\begin{thm}\label{thm:pressure-B}
Let $f$ be a non-exceptional transcendental meromorphic map in the class $\B$, such that $J(f) \setminus \overline{\PP(f)} \neq \emptyset$. Then for every $z \in J(f) \setminus \overline{\PP(f)}$,
\[
\overline{P}(t,z) = \underline{P}(t,z) = P_{hyp}(t) \quad \text{for every }  t  > 0.
\]
\end{thm}
\begin{proof} First we prove $P_{hyp}(t) \leq \underline{P}(t,z)$ for every $z \in J(f) \setminus \overline{\PP(f)}$. Take a transitive isolated conformal repeller $X \subset J(f)$ and a point $z \in J(f) \setminus \overline{\PP(f)}$. 
Note that for every $z_0 \in f^{-1}(z)$ we have $S_n(t,z_0) \leq |f^*(z_0)|^t S_{n+1}(t,z)$, which implies $\underline{P}(t, z_0) \leq \underline{P}(t,z)$, so by induction,
\[
\underline{P}(t, z_0) \leq \underline{P}(t,z) \qquad\text{for every }z_0 \in \bigcup_{n = 1}^\infty f^{-n}(z).
\]
Since $\bigcup_{n = 1}^\infty f^{-n}(z)$ are dense in $J(f)$, we can take $z_0 \in \bigcup_{n = 1}^\infty f^{-n}(z)$ in a small neighbourhood of a point $z_1 \in X$. Then we show $P(f|_X,t) \leq \underline{P}(t,z_0)$ in the same way as in the proof of Theorem~\ref{thm:pressure-S} and we conclude that $P(f|_X,t) \leq \underline{P}(t,z)$, which shows $P_{hyp}(t) \leq \underline{P}(t,z)$.

Now we show $P_{hyp}(t) \geq \overline{P}(t,z)$ for every $z \in J(f) \setminus \overline{\PP(f)}$. First we prove an analogue of Proposition~\ref{prop:S^K}, i.e.~we show that for every $z \in J(f) \setminus \overline{\PP(f)}$,
\begin{equation}\label{eq:S^K-B}
P_{hyp}(t) \geq \sup_K\limsup_{n\to\infty}\frac{1}{n} \ln S^K_n(t,z),
\end{equation}
where the supremum is taken over all compact sets $K \subset J(f)$, which do not contain the exceptional values of $f$. To show \eqref{eq:S^K-B}, take $z_0 \in J(f) \setminus \overline{\PP(f)}$ and 
\[
D = \DD(z_0, \delta),
\]
where $\delta > 0$ is so small that $D \subset J(f) \setminus \overline{\PP(f)}$. Then we define the number $m_n$ and the branches $g_z$, $\tilde g_{w_z}$ in the same way as in the proof of Proposition~\ref{prop:S^K}, replacing $f^{l_n}(z_0)$ by $z_0$. In this way we construct a conformal iterated function system 
\[
\F_n = \{h_z\}_{z\in f^{-n}(z_0)\cap K},
\]
where $h_z = \tilde g_{w_z} \circ g_z$ are suitable inverse branches of $f^{-N}$ for $N = n+m_m$. Taking $\Lambda_n$ to be the $f^N$-invariant conformal repeller defined by $\F_n$ and $Y_n = \bigcup_{j=0}^{N-1} f^j(\Lambda_n)$ to be the suitable $f$-invariant conformal repeller, we show (in the same way as in the proof of Proposition~\ref{prop:S^K}) that $P(f|_{Y_n}, t) \geq \limsup_{n\to\infty}\frac{1}{n} \ln S^K_n(t,z_0)$, concluding the proof of \eqref{eq:S^K-B}.

Having \eqref{eq:S^K-B}, we repeat the proof of Theorem~\ref{thm:pressure-S}, split into three cases. The only difference is that to apply
Corollary~\ref{cor:log} for an exceptional value $a \in J(f)$ we use the assumption that $f$ is non-exceptional (for $f \in \SSS$ it was satisfied automatically).
\end{proof}

By Theorem~\ref{thm:pressure-B}, for non-exceptional $f \in \B$ with $J(f) \setminus \overline{\PP(f)} \neq \emptyset$ we can define the pressure of $f$ as
\[
P(t) = P(f, t, z) = \lim_{n\to\infty} \frac{1}{n} \ln S_n(t, z)
\]
for any point $z \in J(f) \setminus \overline{\PP(f)}$.

\begin{prop}\label{prop:pressure-B}
If $f \in \B$ is non-exceptional and $J(f) \setminus \overline{\PP(f)} \neq \emptyset$, then 
\begin{itemize}
\item 
$P(t) > -\infty$ for every $t > 0$.
\item 
the function $t \mapsto P(t)$ is non-increasing for $t \in (0, +\infty)$
\item
the function $t \mapsto P(t)$ is convex $($and hence continuous$)$ for $t \in (t_0, +\infty)$, where $t_0 = \inf\{t > 0: P(t) \text{ is finite}\}$,
\item  $P(2) \leq 0$.
\end{itemize}
\end{prop}
\begin{proof}
The arguments are the same as in the proof of Proposition~\ref{prop:pressure}.
\end{proof}

In the same way as for functions from $\SSS$, we obtain the following.
\begin{bowen-B}
If $f$ is a non-exceptional transcendental meromorphic map in the class $\B$, such that $J(f) \setminus \overline{\PP(f)} \neq \emptyset$, then 
\[
\dim_H J_r(f)  = \dim_{hyp} J(f) = \delta(f),
\]
where $\delta(f) = \inf \{t > 0: P(t) \leq 0\}$.
\end{bowen-B}

We end this section by looking at the special case of hyperbolic maps
in the class $\B$ (see also \cite{Z}). Then we can prove the existence of the pressure in a more direct way.

\begin{thm}\label{thm:pressure-B-hyp} Let $f$ be a hyperbolic map in the class $\B$. Then for every $t > 0$, we have $P(t) = P(t, z) > -\infty$ for every $z \in J(f)$ $($and the limit in the definition of $P(t, z)$ exists$)$. Moreover, $P(t) > 0$ for every $0 < t < \delta(f)$ and $P(t) < 0$ for every $t > \delta(f)$.
\end{thm}
\begin{proof} Obviously, if $f \in \B$ is hyperbolic, then $f$ is non-exceptional and $J(f) \cap \overline{\PP(f)} = \emptyset$. This together with Theorem~\ref{thm:pressure-B} shows that $P(t) = P(t, z)$ for every $z \in J(f)$. Note that for $z \in J(f)$ we have
\[
S_{n+m}(t,z) \geq c S_n(t,z) S_m(t,z)
\]
for some constant $c > 0$ and every $n,m > 0$, which follows easily from Corollary~\ref{cor:sta}, the Spherical Koebe Distortion Theorem and the hyperbolicity of $f$. Then we conclude that the limit $P(t) = \lim_{n\to\infty} \frac{1}{n} \ln S_n(t,z)$ exists and $P(t) > -\infty$.
The second assertion of the theorem follows immediately from the proofs of \cite[Lemmas~3.3 and~3.5]{S} by G.~Stallard. 
\end{proof}

\end{document}